\documentclass[12pt,a4paper,reqno]{amsart}

\usepackage[active]{srcltx}
\usepackage{enumerate}
\usepackage{float}
\usepackage{amssymb}

\usepackage[all]{xy}

\usepackage{cases}

\newcommand{\Z}{{\mathbb Z}}

\newcommand{\F}{{\mathbb F}}

\renewcommand{\P}[1]{\mathcal{P}^{#1}}

\renewcommand{\ker}{\operatorname{Ker}\nolimits}
\newcommand{\im}{\operatorname{Im}\nolimits}

\newcommand{\A}{\ifmmode{\mathcal{A}}\else${\mathcal{A}}$\fi}
\newcommand{\K}{\ifmmode{\mathcal{K}}\else${\mathcal{K}}$\fi}
\newcommand{\U}{\ifmmode{\mathcal{U}}\else${\mathcal{U}}$\fi}
\newcommand{\M}{\ifmmode{\mathcal{M}}\else${\mathcal{M}}$\fi}
\newcommand{\N}{\ifmmode{\mathcal{N}}\else${\mathcal{M}}$\fi}
\newcommand{\Ff}{\ifmmode{\mathcal{F}}\else${\mathcal{F}}$\fi}
\newcommand{\Ll}{\ifmmode{\mathcal{L}}\else${\mathcal{L}}$\fi}
\newtheorem{theorem}{Theorem}[section]
\newtheorem{proposition}[theorem]{Proposition}
\newtheorem{corollary}[theorem]{Corollary}

\theoremstyle{definition}
\newtheorem{definition}[theorem]{Definition}
\newtheorem{remark}[theorem]{Remark}

\theoremstyle{remark}

\theoremstyle{plain}

\title[Cohomological uniqueness of some $p$-groups]
{Cohomological uniqueness of some $p$-groups}

\author{Antonio D{\'\i}az}
\address{Department of Mathematical Sciences, University of Copenhagen, Universitetsparken 5,
DK-2100 Copenhagen \o,Denmark}
\email{adiaz@math.ku.dk}

\author{Albert Ruiz}
\address{Departament de Matem{\`a}tiques,
Universitat Aut{\`o}noma de Barcelona, 08193 Cerdanyola del
Vall{\`e}s, Spain.}
\email{Albert.Ruiz@uab.cat}

\author{Antonio Viruel}
\address{Departamento de {\'A}lgebra, Geometr{\'\i}a y Topolog{\'\i}a,
Universidad de M{\'a}\-la\-ga, Apdo correos 59, 29080 M{\'a}laga,
Spain.}
\email{viruel@agt.cie.uma.es}

\thanks{
\indent First and third authors are partially supported
by MEC grant MTM2007-60016, Junta de Andaluc{\'\i}a grant FQM-213
and P07-FQM-2863.\\ \indent Second author is partially supported
by MEC grant MTM2007-61545.\\ \indent Second and third authors are
partially supported by Generalitat de Catalunya grant
2009SGR-1092.}

\begin{document}

\begin{abstract}
In this paper we consider classifying spaces of a family of
$p$-groups and we prove that mod~$p$ cohomology enriched with
Bockstein spectral sequences determines their homotopy type among
$p$-completed CW-complexes.
\end{abstract}

\maketitle


\section{Introduction}

Let $p$ be a prime number. A naive way of describing Bousfield-Kan
$p$-completion functor \cite{BK} is to say that it does transform
mod~$p$ cohomology isomorphisms into actual homotopy equivalences.
It is then therefore natural to think that the homotopy type of a
$p$-complete space $X$ should be characterized in some sense by
its mod~$p$ cohomology ring $H^*X$. Classifying spaces of finite
$p$-groups provide nice examples of $p$-complete spaces. Then the
following question arises: given a finite $p$-group $P$, and a
$p$-complete space $X$ such that $H^*X\cong H^*BP$, is then
$X\simeq BP$?

One would like to give a positive answer to the question above,
but the very first step towards that positive answer is to
understand, or to give the appropriate meaning to, the isomorphism
$H^*X\cong H^*BP$.

It is well known that there are infinitely many examples of non
isomorphic finite $p$-groups (hence infinitely many examples of
non homotopic $p$-complete spaces) having isomorphic mod~$p$
cohomology rings, even as unstable algebras (see \cite{CC} for a
general proof of this fact in the case of $p=2$). This is not
surprising since $p$-completion does not invert abstract mod~$p$
cohomology isomorphisms, but just those which are induced by
continuous maps, and these compare unstable algebras plus
secondary operations.

In this direction, Broto-Levi \cite{BL} suggested that mod~$p$
cohomology rings of finite $p$-groups should be considered objects
in the category $\K_\beta$ of unstable algebras endowed with
Bockstein spectral sequences (see Section \ref{Notation} for
precise definitions). Here we follow that line and consider the
family of groups studied by I. Leary in \cite{Leary}, proving:

\begin{theorem}\label{th-main}
Let $p$ be an odd prime and define the finite $p$-group
$$P(p,n)=\langle A,B,C \mid A^p=B^p=C^{p^{n-2}}=[A,C]=[B,C]=1,
[A,B]=C^{p^{n-3}}\rangle.$$ Given $X$ a $p$-complete CW-complex:
\begin{enumerate}[\rm (a)]
\item If $n=3,4$ and $H^*(X)\cong H^*(BP(p,n))$ as unstable algebras, then
$X\simeq BP(p,n)$.
\item If $n\geq 5$ and $H^*_\beta(X)\cong H^*_\beta(BP(p,n))$ as objects in $\K_\beta$, then $X\simeq BP(p,n)$.
\end{enumerate}
\end{theorem}
\begin{proof}
Statement (a) is proved in Corollary \ref{cor:P(p,3)}, $n=3$, and
Corollary \ref{cor:P(p,4)}.(a), $n=4$. Statement (b) is proved in
Corollary \ref{cor:P(p,4)}.(b).
\end{proof}

Besides of its own topological interest, the result above and the
techniques developed in its proof may be appealing from a group
theoretical point of view. First, since the classifying space of a
finite $p$-group is a $p$-complete CW-complex, Theorem
\ref{th-main} provides a cohomological characterization of
$P(p,n)$:

\begin{theorem}\label{th-second}
Let $p$ be an odd prime and $G$ be a finite $p$-group. Then
$G\cong P(p,n)$ if and only if $H^*_\beta(BG)\cong
H^*_\beta(BP(p,n))$.
\end{theorem}

Second, the ideas in the proof of Theorem \ref{th-main} can be
used to obtain a cohomological characterization of $P(p,n)$ as a
complement of some $N\unlhd G$. This characterization can be seen
as a generalization of Tate's cohomological criteria of
$p$-nilpotency \cite{Tate}:

\begin{theorem}\label{th-third}
Let $p$ be an odd prime and $P(p,n)\leq G$ be a group. Then
$P(p,n)$ is a complement for some $N\unlhd G$ if and only if one
of the following holds:
\begin{enumerate}[\rm (a)]
\item $n=3,4$ and there
exists $\psi:H^*(BP(p,n))\to H^*(BG)$ as unstable algebras such
that $(\operatorname{res}\circ\,\psi)|_{H^1_\beta(BP(p,n))}$ is
the identity.
\item $n\geq 5$ there
exists $\psi:H^*_\beta(BP(p,n))\to H^*_\beta(BG)$ in $\K_\beta$
such that $(\operatorname{res}\circ\,\psi)|_{H^1_\beta(BP(p,n))}$
is the identity.
\end{enumerate}
\end{theorem}
\begin{proof}
If $P(p,n)$ is a complement for some $N\unlhd G$, then the induced
projection $G\stackrel{\pi}{\rightarrow} G/N\cong P(p,n)$ gives
rise to a map between classifying spaces
$BG\stackrel{B\pi}{\rightarrow} BP(p,n)$ that provides the desired
cohomological morphism $\psi=B\pi^*$.

The converse in proved in in Proposition \ref{prop:P(p,3)}, in
case $n=3$, and Proposition \ref{prop:P(p,4)} in case $n>3$.
\end{proof}

\noindent\textbf{Organization of the paper:} In Section
\ref{Notation} we introduce the notation we use along paper. In
Section \ref{P(p,n)} the group $P(p,n)$ is defined and the mod~$p$
cohomology ring of its classifying space is described. Finally, in
Section \ref{Main}, we explore endomorphisms of the mod~$p$
cohomology ring of $BP(p,n)$, and we conclude that mod~$p$
cohomology determines the homotopy type of $BP(p,n)$. Finally, in
Section \ref{Groups} we apply the ideas developed in the previous
section to the group theoretical framework.


\section{Definitions and notation} \label{Notation}
We follow the notation and conventions in \cite[Section 2]{BL}. As
our study is done for a fixed odd prime $p$, we just recall the
definitions in this case.

All the spaces considered here have the homotopy type of a
$p$-complete $CW$-complex. Unless otherwise stated $H^*(X)$ refers
to the cohomology of the space $X$ with trivial coefficients in
$\F_p$.

\begin{definition}
Let $p$ be an odd prime and $K$ be an unstable algebra. A
\emph{Bss (Bockstein spectral sequence) for $K$} is a spectral
sequence of differential graded algebras
$\{E_i(K),\beta_i\}_{i=1}^\infty$ where the differentials have
degree one and such that
\begin{enumerate}[(a)]
\item $E_1(K)=K$ and $\beta_1=\beta$ is the primary Bockstein operator.
\item If $x \in E_i(K)^{\mbox{even}}$ and $x^p\neq 0$ in
$E_{i+1}(K)$, $i\geq 1$, then $\beta_{i+1}(x^p)=x^{p-1}\beta_i(x)$.
\end{enumerate}
\end{definition}
We work in the category $\K_\beta$, whose objects are pairs
$(K;\{E_i(K),\beta_i\}_{i=1}^\infty\})$, where $K$ is an unstable
algebra and $\{E_i(K);\beta_i\}_{i=1}^\infty$ is a Bss for $K$. A
morphism $f\colon K\to K'$ in $\K_\beta$ is a family of morphisms
$\{f_i\}_{i=1}^\infty$, where $f_1\colon K\to K'$ is a morphism of
$\A_p$-algebras and for each $i\geq 2$, $f_i\colon E_i(K)\to
E_i(K')$ is a morphism of differential graded algebras, which as
morphism of graded algebras, is induced by $f_{i-1}$.

The mod~$p$ cohomology of a space $X$ is an object of $\K_\beta$
that is denoted by $H_\beta^*(X)$.

\begin{definition}\label{def:determined}
We say that two spaces $X$ and $Y$ are \emph{comparable} if
$H^*_\beta(X)$ and $H^*_\beta(Y)$ are isomorphic objects in the
category $\K_\beta$. We say that $X$ is \emph{determined by cohomology}
if given a space $Y$ comparable to $X$, there is a homotopy
equivalence $X\simeq Y$.
\end{definition}

\begin{definition}\label{def:weakly}
Let $K_\beta$ be an object in $\K_\beta$. Let $K$ be the
underlying unstable algebra over $\A_p$. We say that $K_\beta$ is
\emph{weakly generated by $x_1, \dots, x_n$} if any endomorphism $f$ of
$K_\beta$ such that the restriction of $f$ to the vector subspace
of $K$ generated by $x_1, \dots , x_n$ is an isomorphism, is an
isomorphism in $\K_\beta$.
\end{definition}


\section{The cohomology of some $p$-groups}\label{P(p,n)}

In this section, the $p$-group $P(p,n)$, $p$ an odd prime, and
$n\geq3$, is introduced, and in what follows the notation in
\cite{Leary} is used.

The group
\begin{equation}\label{defP(p,n)}
P(p,n)=\langle A,B,C \mid A^p=B^p=C^{p^{n-2}}=[A,C]=[B,C]=1,
[A,B]=C^{p^{n-3}}\rangle \,
\end{equation}
have order order $p^n$ and fits in a central extension:
\begin{equation}
0\rightarrow \Z/{p^{n-2}} \rightarrow P(p,n) \rightarrow \Z/p \times \Z/p
\rightarrow 0 \,. \label{extension}
\end{equation}

The cohomology of $P(p,n)$ is calculated in \cite{Leary}:
\begin{theorem}[{\cite[Proposition 3, Theorem 7, Proposition 8]{Leary}}]\label{cohoP(3,3)}
$H^*(BP(3,3))$ is generated by elements $y$, $y'$, $x$, $x'$, $Y$,
$Y'$, $X$, $X'$, $z$ with $$ \deg(y)=\deg(y')=1,
\deg(x)=\deg(x')=\deg(Y)=\deg(Y')=2,$$ $$ \deg(X)=\deg(X')=3
\mbox{ and } \deg(z)=6$$ subject to the following relations: $$
yy'=0,xy'=x'y,yY=y'Y'=xy',yY'=y'Y, $$ $$ YY'=xx',Y^2=xY',Y'^2=x'Y,
$$ $$ yX=xY-xx',y'X'=x'Y'-xx', $$ $$ Xy'=x'Y-xY', X'y=xY'-x'Y,$$
$$ XY=x'X,X'Y'=xX',XY'=-X'Y, xX'=-x'X, $$ $$ XX'=0,
x(xY'+x'Y)=-xx'^2, x'(xY'+x'Y)=-x'x^2, $$ $$ x^3y'-x'^3y=0,
x^3x'-x'^3x=0,$$ $$ x^3Y'+x'^3Y=-x^2x'^2 \mbox{ and }
x^3X'+x'^3X=0. $$ Moreover, the action of the mod~$3$ Steenrod
algebra is determined by: $$ \beta(y)=x, \beta(y')=x', \beta(Y)=X,
\beta(Y')=X', $$ $$ \P{1}(X)=x^2X+zy,  \, \P{1}(X')=x'^2X'-zy', \,
\mbox{ and } \, \P{1}(z)=zc_2,$$ where $c_2=xY'-x'Y-x^2-x'^2$.
\end{theorem}

\begin{theorem}[{\cite[Proposition 3, Theorem 6, Proposition 8]{Leary}}]\label{cohoP(p,3)}
For an odd prime $p\geq 5$, the cohomology $H^*(BP(p,3))$ is
generated by elements $y$, $y'$, $x$, $x'$, $Y$, $Y'$, $X$, $X'$,
$d_4$, \dots, $d_p$, $c_4$, \dots, $c_{p-1}$ and $z$ with $$
\deg(y)=\deg(y')=1, \deg(x)=\deg(x')=\deg(Y)=\deg(Y')=2,$$ $$
\deg(X)=\deg(X')=3, \deg(d_i)=2i-1, \deg(c_i)=2i, \mbox{ and }
\deg(z)=2p$$ subject to the following relations: $$ yy'=0,\,
xy'=x'y,\, yY=y'Y'=0,\, yY'=y'Y, $$ $$ Y^2=Y'^2=YY'=0,\, yX=xY,\,
y'X'=x'Y',$$ $$Xy'=2xY'+x'Y,\, X'y=2x'Y+xY'$$ $$ XY=X'Y'=0,\,
XY'=-X'Y,\, xX'=-x'X,$$ $$x(xY'+x'Y)=x'(xY'+x'Y)=0,$$
$$x^py'-x'^py=0,\, x^px'=x'^px=0,$$ $$x^pY'+x'^pY=0,\,
x^pX'+x'^pX=0,$$ $$c_iy=\begin{cases} 0 \nonumber \\ -x^{p-1}y
\nonumber\end{cases} c_iy'=\begin{cases} 0 & \text{for
$i<p-1$}\nonumber \\ -x'^{p-1}y' & \text{for
$i=p-1$,}\nonumber\end{cases} $$ $$c_ix=\begin{cases} 0 \nonumber
\\ -x^p \nonumber\end{cases} c_ix'=\begin{cases} 0 & \text{for
$i<p-1$}\nonumber \\ -x'^p & \text{for
$i=p-1$,}\nonumber\end{cases} $$ $$c_iY=\begin{cases} 0 \nonumber
\\ -x^{p-1}Y \nonumber\end{cases} c_iY'=\begin{cases} 0 &
\text{for $i<p-1$}\nonumber \\ -x'^{p-1}Y' & \text{for
$i=p-1$,}\nonumber\end{cases} $$ $$c_iX=\begin{cases} 0 \nonumber
\\ -x^{p-1}X \nonumber\end{cases} c_iX'=\begin{cases} 0 &
\text{for $i<p-1$}\nonumber \\ -x'^{p-1}X & \text{for
$i=p-1$,}\nonumber\end{cases} $$ $$c_ic_j=\begin{cases} 0 &
\text{for $i+j<2p-2$}\nonumber\\x^{2p-2}+x'^{2p-2}-x^{p-1}x'^{p-1}
& \text{for $i=j=p-1$,}\nonumber\end{cases}$$ $$d_iy=\begin{cases}
0 \nonumber \\ -x^{p-1}Y \nonumber\end{cases} d_iy'=\begin{cases}
0 & \text{for $i<p$}\nonumber \\ -x'^{p-1}Y' & \text{for
$i=p$,}\nonumber\end{cases} $$ $$d_ix=\begin{cases} 0 \nonumber \\
-x^{p-1}y \nonumber \\ x^{p-1}X\nonumber\end{cases}
d_ix'=\begin{cases} 0 & \text{for $i<p-1$}\nonumber \\ -x'^{p-1}y'
& \text{for $i=p-1$}\nonumber \\ -x'^{p-1}X' & \text{for
$i=p$,}\nonumber\end{cases} $$ $$d_iY=0,\, d_iY'=0,$$
$$d_iX=\begin{cases} 0 \nonumber \\ -x^{p-1}Y \nonumber\end{cases}
d_iX'=\begin{cases} 0 & \text{for $i\neq p-1$}\nonumber \\
-x'^{p-1}Y' & \text{for $i=p-1$,}\nonumber\end{cases} $$
$$d_id_j=\begin{cases} 0 & \text{for $i<p$ or
$j<p-1$}\nonumber\\x^{2p-3}Y+x'^{2p-3}Y'+x^{p-1}x'^{p-2}Y' &
\text{for $i=p$ and $j=p-1$,}\nonumber\end{cases}$$
$$d_ic_j=\begin{cases} 0 & \text{for $i<p-1$ or $j<p-1$}\nonumber
\\ -x^{2p-3}y+x'^{2p-3}y'-x^{p-1}x'^{p-2}y' & \text{for
$i=j=p-1$}\nonumber \\ -x^{2p-3}X+x'^{2p-3}X'-x^{p-1}x'^{p-2}X' &
\text{for $i=p$, $j=p-1$.}\nonumber\end{cases} $$ Moreover, the
action of the mod~$p$ Steenrod algebra is determined by: $$
\beta(y)=x,\, \beta(y')=x',\, \beta(Y)=X,\, \beta(Y')=X', $$ $$
\beta(d_i)=\begin{cases} c_i & \text{for $i<p$}\nonumber\\ 0 &
\text{for $i=p$.}\nonumber\end{cases} $$ $$
\P{1}(X)=x^{p-1}X+zy,\, \P{1}(X')=x'^{p-1}X'-zy', $$ $$
\P{1}(c_i)=
\begin{cases} izc_{i-1} & \text{ if $2\leq i<p-1$} \nonumber \\
-zc_{p-2}+x^{2p-2}+x'^{2p-2}-x^{p-1}x'^{p-1} & \text{ if $i=p-1$,}
\end{cases}
$$ $$ \P{1}(z)=zc_{p-1}, $$ where $c_2$ and $c_3$ are non zero
multiples of $xY'+x'Y$ and $XX'$ respectively.

\end{theorem}

\begin{remark}\label{remark:cohoP33}
As stated in \cite[p.\ 71]{Leary} one can verify that in the
cohomology ring $H^*(BP(p,3))$, $p\geq5$, any product of the
generators $y$, $y'$, $x$, $x'$, $Y$, $Y'$, $X$, $X'$ in degree
greater than $6$ may be expressed in the form $$
\begin{array}{ll}
f_1+f_2Y+f_3Y' & \text{for even total degree}\\
f_1y+f_2y'+f_3X+f_4X' & \text{for odd total degree,}\\
\end{array}
$$ where each $f_i$ is a polynomial in $x$ and $x'$. So, for
$n\leq p$, any element $u \in H^{2n-1}(BP(p,3))$ can be expressed
as: $$ u=ad_n+f_1y+f_2y'+f_3X+f_4X' , $$ where $a \in \F_p$ and
each $f_i$ is a polynomial in $x$ and $x'$.
\end{remark}

\begin{theorem}[{\cite[Theorem 4]{Leary}}]\label{cohoP(p,n)}
For $n\geq 4$, $H^*(BP(p,n))$ is generated by elements
$u$, $y$, $y'$, $x$, $x'$, $c_2$, $c_3$,\ldots, $c_{p-1}$, $z$ with
$$
\deg(u)=\deg(y)=\deg(y')=1, \,\deg(x)=\deg(x')=2,\, \deg(c_i)=2i,\, \deg(z)=2p,
$$
subject to the following relations:
$$ xy'=x'y,\,x^py'=x'^py,\,x^px'=x'^px, $$
$$ c_iy=\begin{cases} 0 \nonumber \\ -x^{p-1}y \nonumber\end{cases}
c_iy'=\begin{cases} 0 & \text{for $i<p-1$}\nonumber \\ -x'^{p-1}y' & \text{for $i=p-1$,}\nonumber\end{cases}
$$
$$ c_ix=\begin{cases} 0 \nonumber \\ -x^p \nonumber\end{cases}
c_ix'=\begin{cases} 0 & \text{for $i<p-1$} \nonumber \\ -x'^p & \text{for $i=p-1$,}\nonumber\end{cases}
$$
$$ c_ic_j=\begin{cases} 0 & \text{for $i+j<2p-2$} \nonumber \\x^{2p-2}+x'^{2p-2}-x^{p-1}x'^{p-1}  & \text{for $i=j=p-1$.}\nonumber\end{cases}$$

Moreover, we have the following operations of the mod~$p$ Steenrod
algebra: $$ \beta(y)=x, \beta(y')=x', \beta(u)=\begin{cases} 0 &
\text{for $n>4$} \nonumber \\ y'y & \text{for $n=4$},
\nonumber\end{cases} $$ and $$ \P{1}(z)=zc_{p-1},
\P{1}(c_i)=\begin{cases} izc_{i-1} & \text{for $i<p-1$}
\nonumber\\ -zc_{p-2}+x^{2p-2}+x'^{2p-2}-x^{p-1}x'^{p-1} &
\text{for $i=p-1$,}\nonumber\end{cases} $$ where $c_1=y'y$.
\end{theorem}

\begin{remark}\label{remark:bockstein}
If we look at $H^*_\beta(BP(p,n))$ for $n\geq4$ as an object in
$\K_\beta$ we have that there is a Bockstein operator
$\beta_{n-3}(u)=yy'$ \cite[p.\ 66]{Leary}. From this and the
cohomology of $BP(p,n)$ we can deduce that for $n=3$ and $n=4$ the
study of the cohomological uniqueness of these spaces can be done
in $\K$, the category of unstable algebras, instead of working in
$\K_\beta$.
\end{remark}
\begin{remark}\label{remark:tower}
Consider the groups $B\Z/p^{i}\times B\Z/p\times B\Z/p$ and fix
the following notation for the cohomology: $$ H^*(B\Z/p^{i}\times
B\Z/p\times B\Z/p;\F_p)=E(u_i,y,y')\otimes\F_p[v_i,x,x'] $$ where
generators are sorted as components.

There is a tower of extensions:
\begin{multline} \nonumber
BP(p,n) \stackrel{\pi_{n-3}}{\to} B\Z/p^{n-3}\times B\Z/p\times
B\Z/p \stackrel{\pi_{n-4}}{\to} B\Z/p^{n-4}\times B\Z/p\times
B\Z/p \to \cdots \\ \cdots \stackrel{\pi_{1}}{\to} B\Z/p\times
B\Z/p\times B\Z/p
\end{multline}
where each extension $\pi_i$ for $i<n-3$: $$ 0 \rightarrow B\Z/p
\rightarrow B\Z/{p^i}\times B\Z/p\times B\Z/p
\stackrel{\pi_{i}}{\rightarrow} B\Z/{p^{i-1}}\times B\Z/p\times
B\Z/p \rightarrow 0 $$ is classified by $\beta_{i-1}(u_i)$ and $$
0 \rightarrow B\Z/p\rightarrow BP(p,n)\stackrel{\pi_{n-3}}{\to}
B\Z/p^{n-3}\times B\Z/p\times B\Z/p \rightarrow 0 $$ is classified
by $\beta_{n-3}(u_{n-3})-yy'$.
\end{remark}


\section{Cohomological uniqueness}\label{Main}
Let $p$ be an odd prime, $n\geq 3$ and $P(p,n)$ be the group
defined in Equation \eqref{defP(p,n)}. In this section we prove
that the homotopy type of the classifying space of $P(p,n)$ is
determined by its cohomology (Definition \ref{def:determined}).
For $n\leq 4$ we do not need to use higher Bocksteins and it is
enough to consider the structure of unstable algebra.

\begin{theorem}\label{isoP(3,3)}
Let $\varphi \colon H^*(BP(3,3))\rightarrow H^*(BP(3,3))$ be a
homomorphism of $\A_3$-algebras which restricts to the identity in
$H^1$. Then $\varphi$ is an isomorphism.
\end{theorem}
\begin{proof} In this proof we follow the notation in Theorem \ref{cohoP(3,3)} for generators
and relations in cohomology.

By hypothesis $\varphi(y)=y$ and $\varphi(y')=y'$. Now, since
$\beta(y)=x$ and $\beta(y')=x'$, then
$\varphi(x)=\varphi(\beta(y))=\beta(\varphi(y))=\beta(y)=x$ and
analogously $\varphi(x')=x'$. Moreover, by dimensional reasons, $$
\varphi(Y)=aY+bY'+cx+dx' $$ for some $a,b,c,d \in \F_3$. Because
$yY=xy'$ we obtain $$
xy'=\varphi(xy')=\varphi(yY)=y\varphi(Y)=ayY+byY'+cyx+dyx' $$ and
regrouping terms $$xy'=(a+d)xy'+byY'+cyx . $$ From here we obtain
$a+d=1$ and $b=c=0$, and $\varphi(Y)=aY+dx'$ with $a+d=1$.
Analogously $\varphi(Y')=bY'+cx$ with  $b,c \in \F_3$ and $b+c=1$.
Now, as $Y^2=xY'$, we have
\begin{align}
\varphi(Y)^2&=x\varphi(Y') \notag \\
a^2Y^2+d^2x'^2+2adYx'&=bxY'+cx^2 .\notag
\end{align}
This implies that $c=d=0$ and $a^2=a=b=1$. So $\varphi(Y)=Y$ and
$\varphi(Y')=Y'$, and applying Bockstein again $\varphi(X)=X$ and
$\varphi(X')=X'$ too. So $\varphi$ is the identity up to dimension
five and it remains to check where does it map $z$.

Using the first Steenrod power of $X$
\begin{align}
\varphi(\P{1}(X))&=\P{1}(\varphi(X)) \notag \\ \varphi(x^2X+zy)
&=\P{1}(X) \notag \\ x^2X+\varphi(z)y &=x^2X+zy \notag \\
\varphi(z)y      &=zy.    \notag
\end{align}
Thus $\varphi(z)=z+\alpha$ where $\alpha y =0$ and
$\alpha\in\langle y,y',x,x',Y,Y',X,X' \rangle$. So
$\varphi(\alpha)=\alpha$, $z=\varphi(z-\alpha)$ and $\varphi$ is
an epimorphism. In fact, because $H^*(BP(3,3))$ is a finite
dimensional $\F_3$-vector space in each dimension, $\varphi$ is an
isomorphism dimension-wise, and thus $\varphi$ is an isomorphism.
\end{proof}

\begin{theorem}\label{isoP(p,3)}
Let $p\geq 5$ be a prime. If $\varphi:H^*(BP(p,3))\rightarrow
H^*(BP(p,3))$ is a homomorphism in $\K$ that restricts to the
identity in $H^1$, then $\varphi$ is an isomorphism.
\end{theorem}
\begin{proof}
Consider the notation of generators and relations in
$H^*_\beta(BP(p,3))$ given in Theorem {\rm\ref{cohoP(p,3)}}. We
calculate the image under $\varphi$ of every generator in
$H^*(BP(p,3))$.

As $\varphi$ is the identity on $y$ and $y'$, applying Bockstein
operations we get that $\varphi(x)=x$ and $\varphi(x')=x'$.

As $Y$ is in degree $2$, there exist coefficients $a$, $b$, $c$,
$d$ such that $$ \varphi(Y)=ax+bx'+cY+dY'.$$ Using the relation
$Y^2=0$, we get $\varphi(Y)^2=0$, which implies that $a=b=0$, and
so $\varphi(Y)=cY+dY'$. The relation $yY=0$ implies
$0=y\varphi(Y)=dyY'$, so $d=0$, getting that there is $c$ such
that $\varphi(Y)=cY$. Using the same arguments, there is $d$ such
that $\varphi(Y')=dY'$.

According to Remark \ref{remark:cohoP33}, there are $a_n \in \F_p$
and $f_{n,i}$ polynomials in $x$ and $x'$ such that for $4\leq
n\leq p$ $$
\varphi(d_n)=a_nd_n+f_{n,1}y+f_{n,2}y'+f_{n,3}X+f_{n,4}X' , $$ and
applying the Bockstein operation, we get that for $4\leq n\leq
p-1$: $$ \varphi(c_n)=a_nc_n+f_{n,1}x+f_{n,2}x'. $$

The relation $c_{p-1}x=-x^p$ gives rise to the following equalities:
$$
\begin{array}{ll}
-x^p & =\varphi(-x^p) =\varphi(c_{p-1}x)=\varphi(c_{p-1})\varphi(x)=\varphi(c_{p-1})x= \\
  &  =
  a_{p-1}c_{p-1}x+f_{p-1,1}x^2+f_{p-1,2}xx'=-a_{p-1}x^p+f_{p-1,1}x^2+f_{p-1,2}xx',
\end{array}
$$ so $(a_{p-1}-1)x^{p}=f_{p-1,1} x^2+f_{p-1,2} x x'$, and as
there are no relations involving $x$ and $x'$ till degree $2p+2$,
we can simplify:
\begin{equation}\label{eq:f1xf2xprima}
f_{p-1,1}x+f_{p-1,2}x'=(a_{p-1}-1)x^{p-1}.
\end{equation}

Doing the same computations using the relation $c_{p-1}x'=-x'^p$ we get
\begin{equation}\label{eq:f1xf2xprimaprima}
f_{p-1,1}x+f_{p-1,2}x'=(a_{p-1}-1)x'^{p-1}.
\end{equation}
Comparing now \eqref{eq:f1xf2xprima} and
\eqref{eq:f1xf2xprimaprima}, and using again that there is no
relation between $x$ and $x'$ till degree $2p+2$, we get
$a_{p-1}=1$, $\varphi(c_{p-1})=c_{p-1}$.

Now we see that $\varphi(c_n)=a_nc_n$, for $4\leq n<p-1$: using
the relation $c_nx=0$ and applying $\varphi$ we get
$f_{n,1}x+f_{n,2}x'=0$, so
\begin{equation}\label{eq:fcn}
\varphi(c_n)=a_nc_n .
\end{equation}

In order to calculate $\varphi(z)$, we apply $\varphi$ to the
equality: $$
\P{1}(c_{p-1})=-zc_{p-2}+x^{2p-2}+x'^{2p-2}-x^{p-1}x'^{p-1}.$$
Since $\varphi(c_{p-1})=c_{p-1}$, $\varphi(x)=x$, and
$\varphi(x')=x'$, we get
\begin{equation}\label{eq:zcpmenos2}
zc_{p-2}=\varphi(z)a_{p-2}c_{p-2} .
\end{equation}
As the generator $z$ does not appear in any relation, Equation
\eqref{eq:zcpmenos2} implies that $a_{p-2}\neq0$ and
$\varphi(z)=a_{p-2}^{-1}z+g$, where $g$ is an expression not
involving $z$, and such that $gc_{p-2}=0$.

We use that we know that $a_{p-2}\neq0$ to check that $a_n\neq 0$
for $4\leq n<p-2$ with an induction argument: assume
$\varphi(c_n)=a_nc_n$ with $a_n\neq0$ and $5\leq n \leq p-2$, and
compute $\varphi(c_{n-1})$: $$
nzc_{n-1}=\P{1}(c_n)=\P{1}(\varphi(a_n^{-1}c_n))=\varphi(a_n^{-1}\P{1}(c_n))=
a_n^{-1}n\varphi(z)a_{n-1}c_{n-1} . $$ This implies
$zc_{n-1}=a_n^{-1}a_{n-1}\varphi(z)c_{n-1}$, and this can only
happens if $a_{n-1}\neq 0$ and $\varphi(z)=a_{n}a_{n-1}^{-1}z+g$
($g$ not involving $z$).

From the expression $c_3=\mu XX'$ we deduce that
$\varphi(c_3)=a_3c_3$ with $a_3=cd$, where $c$ and $d$ were
introduced at the beginning of the proof and are such that
$\varphi(Y)=cY$ and $\varphi(Y')=dY'$. The argument above for
$\P1(c_4)$ shows that $a_3$ is non-zero neither. Hence $c$, $d$
and $a_n$ for all $n\in\{3,\ldots,p-1\}$ are non-zero.

Let us check now that the coefficients $c$ and $d$ are equal: recall that $c_2$ was defined as
$\lambda (xY'+x'Y)$ with $\lambda$ non-zero. Then, applying $\P{1}$ to $c_3$ we get:
$$
\begin{array}{ll}
3zc_2 & =\P{1}(c_3)=\P{1}(\varphi(a_3^{-1}c_3))=a_3^{-1}\varphi(\P{1}(c_3))=a_3^{-1}(\varphi(3zc_2)) =\\
 & = a_3^{-1}3\varphi(z)\varphi(c_2) =
 a_3^{-1}3\varphi(z)\lambda(dxY'+cx'Y),
\end{array}
$$
which implies $\lambda z(xY'+x'Y)=a_3^{-1}\lambda\varphi(z)(dxY'+cx'Y)$ and can be simplified to:
\begin{equation}\label{eq:zxYprima}
 zxY'+zx'Y=da_3^{-1}\varphi(z)xY'+ca_3^{-1}\varphi(z)x'Y .
\end{equation}
Again, as $z$ does not appear in any relation, Equation \eqref{eq:zxYprima} can be true only if $c=d$.
In particular, $\varphi(c_2)=a_2c_2$ with $a_2=c\lambda\neq 0$.

Now we can assume that all the coefficients $a_n$ for $2\leq n\leq
p-1$ and  $c$ and $d$ are equal to $1$: as all are different to
zero, and $r^{p-1}=1$ if $r\in \F_p\setminus\{0\}$,
$\varphi^{p-1}$ is the identity in $Y$, $Y'$ and $c_n$. Use now
that $\varphi$ is an isomorphism if and only if $\varphi^{p-1}$ is
so. Therefore at this point we have that: $$ \varphi(y)=y,
\varphi(y')=y',\, \varphi(x)=x,\, \varphi(x')=x',$$ $$
\varphi(Y)=Y,\, \varphi(Y')=Y',\, \varphi(X)=X,\, \varphi(X')=X',
$$ $$ \varphi(c_i)=c_i \text{ for $2\leq i\leq p-1$ },\,
\varphi(d_i)=d_i+g_i \text{ for $4\leq i\leq p-1$, and } f(z)=z+g
$$ where $g$ and all $g_i$ are expressions in $x$, $x'$, $y$,
$y'$, $X$, $X'$, $Y$ and $Y'$. This implies that all generators
but $d_p$ are in the image of $\varphi$.

The image of $d_p$, as it is in odd degree greater than 6, must
be: $$ \varphi(d_p)=a_p d_p + f_{p,1}y +
f_{p,2}y'+f_{p,3}X+f_{p,4}X' $$ with $a_p \in \F_p$, and $f_{p,i}$
polynomials in $x$ and $x'$. As $\beta(d_p)=0$, the Bockstein
operation on $\varphi(d_p)$ must vanish, and this means: $$
0=\beta(\varphi(d_p))=f_{p,1}x+f_{p,2}x'.$$ So this is a
polynomial in $x$, $x'$ which must be zero. As there are not
relations involving just $x$ and $x'$ in this degree, we deduce
that there exist $f_p$ a polynomial in $x$ and $x'$ such that
$f_{p,1}=f_p x'$ and $f_{p,2}=f_p x$. This implies that (recall
$xy'=x'y$), $$ f_{p,1}y + f_{p,2}y'=f_p(x'y-xy')=0,$$ and then $$
\varphi(d_p)=a_p d_p + f_{p,3}X+f_{p,4}X'.$$ As any expression on
$x$, $x'$, $X$ and $X'$ is in the image, we have only to check
that $a_p\neq 0$. To do that we assume that $a_p=0$, getting a
contradiction.

If $\varphi(d_p)=f_{p,3}X+f_{p,4}X'$, using the relation $d_p Y=0$ we get
$(f_{p,3}X+f_{p,4}X')Y=0$. Applying that $XY=0$ this means that $f_{p,4}X'Y=0$.
Recall now that $f_{p,4}$ is a polynomial on $x$, $x'$ of degree $2p-4$.
Looking at the relations under this degree, $f_{p,4}X'Y=0$ means that $f_{p,4}$ must be zero.
The same argument can be applied to the relation $d_p Y'=0$, obtaining that
$f_{p,3}=0$. So at this point we have that
$\varphi(d_p)=0$.

Use now that $d_p x = x^{p-1}X \neq 0$, and applying $\varphi$ to
both sides of this equality we get the following contradiction: $$
0 = \varphi(d_p x)= \varphi(x^{p-1}X)= x^{p-1}X \neq 0 . $$
\end{proof}

\begin{theorem}\label{isoP(p,n)}
Let $p$ be an odd prime and consider the notation of the
generators and relations in $H^*_\beta(BP(p,n))$ as in Theorem
{\rm\ref{cohoP(p,n)}}.
\begin{enumerate}[\rm (a)]
\item If $\varphi\colon H^*(BP(p,4))\rightarrow H^*(BP(p,4))$ is a
homomorphism of unstable algebras that fixes $y$ and $y'$, then
$\varphi$ is an isomorphism.
\item If $n\geq 5$ and $\varphi\colon H^*_\beta(BP(p,n))\rightarrow H^*_\beta(BP(p,n))$
is a homomorphism in $\K_\beta$ which fixes $y$
and $y'$. Then $\varphi$ is an isomorphism.
\end{enumerate}
\end{theorem}
\begin{proof}
We prove both results at the same time. Just observe that the
Bockstein used in the proof is $\beta_{n-3}$, which is part of the
mod~$p$ Steenrod algebra when $n=4$.

Starting from $\varphi(y)=y$ and $\varphi(y')=y'$ and using the
Bockstein operator we reach $\varphi(x)=x$ and $\varphi(x')=x'$.
On the other hand there exist $a,b,c\in \F_p$ such that
$\varphi(u)=au+by+cy'$. From Remark \ref{remark:bockstein} we
deduce that $\beta_{n-3}(u)= y'y$, so, as the morphism is in
$\K_\beta$, $$ \varphi(\beta_{n-3}(u))=\beta_{n-3}(\varphi(u))
\Rightarrow y'y = \begin{cases} a y' y + bx +cx' & \text{for
$n=4$}\nonumber \\ a y'y & \text{for $n>4$.}\nonumber\end{cases}
$$ We obtain then that $a=1$ and $b=c=0$ for $n=4$, and that $a=1$
for $n>4$. Hence, $u=\varphi(u)-by-cy'=\varphi(u-by-cy')$ and
$\langle u,y,y',x,x'\rangle \leq \im \varphi$.

Now consider the generator $c_{p-1}$. We can write
$\varphi(c_{p-1})=a_{p-1}c_{p-1}+b x^{p-1}+g_{p-1}$ with $a_{p-1},b\in \F_p$
and $g_{p-1}$ not containing multiples of the monomials $c_{p-1}$ and $x^{p-1}$.
Applying $\varphi$ to the equation  $c_{p-1}x'=-x'^p$ we obtain
$-x'^p=a_{p-1}c_{p-1}x'+b x^{p-1}x'+g_{p-1}x'=-a_{p-1}x'^p+b x^{p-1}x'+g_{p-1}x'$.
The only equation at degree $p$ involving $x'^p$ is $c_{p-1}x'=-x'^p$
(notice that a multiple of the equation $xy'=x'y$ in degree $3$ does not involve $x'^p$ and that this may
occur just for $p=3)$. As $g_{p-1}$ does not contain neither $c_{p-1}$ nor
$x^{p-1}$ we deduce that $a_{p-1}=1$, $b=0$. Hence
$\varphi(c_{p-1})=c_{p-1}+g_{p-1}$.

Next we deal with $c_{p-2}$ of degree $2(p-2)$ and $z$ of degree $2p$. Their images are
$\varphi(c_{p-2})=a_{p-2}c_{p-2}+g_{p-2}$ and
$\varphi(z)=a_z z+g_z$, with $a_{p-2},a_z\in \F_p$, and
$g_{p-2}$ and $g_z$ not involving the monomials $c_{p-2}$ and $z$ respectively.
Write the Steenrod operation $\P{1}(c_{p-1})=-zc_{p-2}+x^{2p-2}+x'^{2p-2}-x^{p-1}x'^{p-1}$ as
$\P{1}(c_{p-1})=-zc_{p-2}+f$, with $f=x^{2p-2}+x'^{2p-2}-x^{p-1}x'^{p-1}$. Applying $\varphi$ we get:
\begin{align}
\varphi(\P{1}(c_{p-1}))&=\P{1}(\varphi(c_{p-1})) \notag \\
\varphi(-zc_{p-2}+)&=\P{1}(c_{p-1}+g_{p-1}) \notag\\
-(a_z z+g_z)(a_{p-2}c_{p-2}+g_{p-2})+f&=-zc_{p-2}+f+\P{1}(g_{p-1})\notag\\
-a_z a_{p-2}zc_{p-2}-a_z zg_{p-2}-a_{p-2}g_zc_{p-2}-g_zg_{p-2}&=-zc_{p-2}+\P{1}(g_{p-1}).\notag
\end{align}
Notice that there is no relation involving the generator $z$ and
the equations involving $c_{p-2}$ are
$c_{p-2}y=c_{p-2}y'=c_{p-2}x=c_{p-2}x'=c_{p-2}c_j=0$ for $j<p$.
Also, the monomial $zc_{p-2}$ cannot appear in $zg_{p-2}$,
$g_zc_{p-2}$, and $g_zg_{p-2}$. Finally, $\P{1}(g_{p-1})$ does not
involve $zc_{p-2}$ as $g_{p-1}$ does not involve $c_{p-1}$ and the
action of $\P{1}$ on $u,y,y',x,x'$ is determined by the axioms.
Hence, $a_z a_{p-2}=1$ and both $a_z$ and $a_{p-2}$ are non-zero.

For the rest of the generators $c_i$ for $i=2,3,\ldots,p-3$ we can write
$\varphi(c_i)=a_ic_i+g_i$, with $a_i\in \F_p$ and $g_i$ not
involving $c_i$. The Steenrod operation $\P{1}(c_{i+1})=(i+1)zc_i$ provides then
\begin{align}
\varphi(\P{1}(c_{i+1}))&=\P{1}(\varphi(c_{i+1})) \notag \\
\varphi((i+1)zc_i)&=\P{1}(\alpha_{i+1}c_{i+1}+g_{i+1}) \notag\\
(i+1)(a_z z+g_z)(a_i c_i+g_i)&=(i+1)a_{i+1}zc_i+\P{1}(g_{i+1})\notag\\
(i+1)(a_z a_i z c_i +a_z z g_i+a_i a_zc_i+a_zg_i)&=(i+1)a_{i+1}zc_i+\P{1}(g_{i+1}).\notag
\end{align}
Notice again that there is no relation involving the generator $z$
and the relations involving $c_i$ are
$c_iy=c_iy'=c_ix=c_ix'=c_ic_j=0$ for $j<2p-2-i$. Also, the
monomial $zc_i$ cannot appear in $zg_i$, $g_zc_i$, and $g_zg_i$.
Moreover, $\P{1}(g_{i+1})$ does not involve $zc_i$ as $g_{i+1}$
does not involve $c_{i+1}$. We deduce that $(i+1)a_z
a_i=(i+1)a_{i+1}$. As $a_z\neq 0$ and $a_{p-2}\neq 0$, an
inductive argument shows that $a_i\neq 0$ for for
$i=2,3,\ldots,p-3$, and hence for all $i=2,3,\ldots,p-1$.

To finish we show that all the generators
$c_2,c_3,\ldots,c_{p-1},z$ are in the image of $\varphi$. We start
with $c_2=\frac{1}{\alpha_2}(\varphi(c_2)-g_2)$. As $g_2\in\langle
u,x,x,y,y'\rangle\leq \im \varphi$ then $c_2$ is also in the image
of $\varphi$. An inductive argument shows that
$c_i=\frac{1}{\alpha_i}(\varphi(c_i)-g_i)$ is in the image of
$\varphi$ as $g_i$ belongs to $\langle
u,x,x',y,y',c_2,c_3,\ldots,c_{i-1}\rangle$. This argument also
applies to show that $z\in \im \varphi$.

Hence, $\varphi$ is an epimorphism. Because $H^*_\beta(BP(p,n))$ is finite in each dimension $\varphi$
is an isomorphism.
\end{proof}

\begin{corollary}
$H^*_\beta(BP(p,n))$ for odd $p$ and $n\geq 3$ is weakly generated
(Definition \rm{\ref{def:weakly}}) by $y$, and $y'$.
\end{corollary}
\begin{proof}
Let $\varphi$ be and endomorphism of $H^*_\beta(BP(p,n))$ which is
an isomorphism on $\langle y,y' \rangle$. Using the outer
automorphism group of $P(p,n)$ which is described in \cite[Lemma
A.5]{DRV} there is a morphism $f \colon BP(p,n) \to BP(p,n)$ such that
the composition $f^* \circ \varphi$ fixes $y$ and $y'$. Use now
Theorems \ref{isoP(3,3)}, \ref{isoP(p,3)} and \ref{isoP(p,n)} to get the result.
\end{proof}

From these results we will obtain the cohomology uniqueness of the
classifying space $BP(p,n)$. We split this result into two
corollaries because the structure of $P(p,3)$ is essentially
different from that of $P(p,n)$, $n>4$.
\begin{corollary}\label{cor:P(p,3)}
Let $p$ be an odd prime and $X$ be a $p$-complete space such that
$H^*(X)\cong H^*(BP(p,3))$ as unstable algebras. Then $X\simeq
BP(p,3)$.
\end{corollary}
\begin{proof}
We begin with the central extension $$ 0\rightarrow \Z/p
\rightarrow P(p,3) \stackrel{\pi} \longrightarrow \Z/p\times\Z/p
\rightarrow 0 $$ classified by $yy'\in H^2(B\Z/p \times B\Z/p)$,
which gives raise to the principal fibration $$ BP(p,3)
\stackrel{B\pi} \longrightarrow B\Z/p\times B\Z/p
\stackrel{yy'}\longrightarrow B^2\Z/p \,. $$ Consider the map
$\pi_X\colon X\rightarrow B\Z/p\times B\Z/p$ that classifies the
classes $y,y'\in H^1(X)$. Then the composite $$ X \stackrel{\pi_X}
\longrightarrow B\Z/p\times B\Z/p \stackrel{yy'} \longrightarrow
B^2\Z/p $$ is null-homotopic because of Theorems \ref{cohoP(3,3)}
and \ref{cohoP(p,3)}, and so $\pi_X$ lifts to $\varphi\colon
X\rightarrow BP(p,3)$, giving the commutative diagram $$
\xymatrix{
  & & BP(p,3) \ar[d]^{B\pi} \\
X \ar@{.>}[rru]^{\varphi} \ar[rr]^{\pi_X} & & B\Z/p\times B\Z/p }
$$ which implies that $\varphi^*$ fixes $y$ and $y'$. Now apply
Theorems \ref{isoP(3,3)} and \ref{isoP(p,3)} to $\varphi^*$.
\end{proof}

\begin{corollary}\label{cor:P(p,4)}
Let $p$ be an odd prime and $X$ be a $p$-complete space.
\begin{enumerate}[\rm (a)]
\item If $H^*(X)\cong H^*(BP(p,4))$ as unstable algebras then
$X\simeq BP(p,4)$.
\item If $n\geq 5$ and $H^*_\beta(X)\cong H^*_\beta(BP(p,n))$ as objects in $\K_\beta$. Then $X\simeq BP(p,n)$.
\end{enumerate}
\end{corollary}
\begin{proof}
Consider the central extension and notations as in Remark
\ref{remark:tower}.

 If $n=4$ we have $$ 0\rightarrow \Z/p
\rightarrow P(p,4) \stackrel{\pi_1}{\rightarrow}
\Z/p\times\Z/p\times\Z/p \rightarrow 0 $$ Let now $\pi_{1,X}$ be
the map $\pi_{1,X}\colon X\rightarrow B\Z/p\times B\Z/p\times
B\Z/p$ that classifies the classes $u,y,y'\in H^1(X)$.

The composite $$ X \stackrel{\pi_{1,X}} \longrightarrow B\Z/p
\times B\Z/p\times B\Z/p \stackrel{yy'-\beta(u)} \longrightarrow
B^2\Z/p $$ is null-homotopic because of Remark \ref{remark:tower},
and so $\pi_{1,X}$ lifts to $\varphi\colon X\rightarrow BP(p,4)$,
giving the commutative diagram $$ \xymatrix{
  & & BP(p,4) \ar[d]^{B\pi_1} \\
X \ar@{.>}[rru]^{\varphi} \ar[rr]^-{\pi_{1,X}} & & B\Z/p\times
B\Z/p\times B\Z/p } $$ This implies that $\varphi^*$ fixes $y$ and
$y'$ and the result is a consequence of Theorem \ref{isoP(p,n)}.

If $n>4$, we must consider before the extensions: $$ 0\rightarrow
\Z/p \longrightarrow \Z/p^{i}\times\Z/p\times\Z/p
\stackrel{\pi_{i-1}}{\longrightarrow}
\Z/p^{i-1}\times\Z/p\times\Z/p \rightarrow 0 , $$ classified by
$\beta_{i-1}(u) \in H^*_\beta(B\Z/p^{i-1}\times B\Z/p\times B
\Z/p)$.

Let $\pi_{1,X}$ be now the map $\pi_{1,X}\colon X\rightarrow
B\Z/p\times B\Z/p\times\Z/p$ that classifies the classes
$u,y,y'\in H^1(X)$. As $\beta_1(u)=0 \in H^*_\beta(X)$ the map
$\pi_{1,X}$ extends to a map $\pi_{2,X}$ having the following
commutative diagram: $$ \xymatrix{
  & & B\Z/p^2\times B\Z/p \times B\Z/p \ar[d]^{B\pi_1} \\
X \ar@{.>}[rru]^{\pi_{2,X}} \ar[rr]^-{\pi_{1,X}} & & B\Z/p\times
B\Z/p\times B\Z/p
} $$ Using the same argument we can proceed extending the map till
$B\pi_{n-3} \colon B\Z/p^{n-3} \times B\Z/p \times B\Z/p$. To do
the last step we use again Remark \ref{remark:tower} and that
$\beta_{n-3}(u)=yy' \in H^*_\beta(X)$, and obtain a map $\varphi$
which gives the commutative diagram: $$ \xymatrix{
  & & BP(p,n) \ar[d]^{B\pi_{n-2} \circ \cdots \circ B\pi_1} \\
X \ar@{.>}[rru]^{\varphi} \ar[rr]^-{\pi_{1,X}} & & B\Z/p\times
B\Z/p\times B\Z/p} $$ which implies again that $\varphi^*$ fixes
$y$ and $y'$, so apply Theorem \ref{isoP(p,n)}.
\end{proof}


\section{Some applications to Group Theory}\label{Groups}

The techniques used along the proof of the last two corollaries in
the previous section can be used to obtain a cohomological
characterization of $P(p,n)$ as a complement of some $N\unlhd G$,
for a super group $P(p,n)\leq G$. Recall that given a group $G$,
and a normal subgroup $N\unlhd G$, $K\leq G$ is a complement for
$N$ if $G=NK$ and $N\cap K=1$, that is, if $G=N\rtimes K$.

Again, we consider the case $n=3$ separately.

\begin{proposition}\label{prop:P(p,3)}
Let $p$ be an odd prime and $G$ be a finite group such that
$P(p,3)\leq G$ and there exists $\psi:H^*(BP(p,3))\to H^*(BG)$ as
unstable algebras such that
$(\operatorname{res}\circ\,\psi)|_{H^1_\beta(BP(p,3))}$ is the
identity. Then $P(p,3)$ is a complement for some $N\unlhd G$.
\end{proposition}
\begin{proof}
As it was announced, we work along the lines in the proof of
Corollary \ref{cor:P(p,3)}. We begin by considering the map
$B\pi_G\colon BG\rightarrow B\Z/p\times B\Z/p$ that classifies the
classes $\psi(y),\psi(y')\in H^1(BG)$. Then
$B\pi_G^*(yy')=B\pi_G^*(y)B\pi_G^*(y')=\psi(y)\psi(y')=\psi(yy')=\psi(0)=0$
(Theorems \ref{cohoP(3,3)} and \ref{cohoP(p,3)}), and the
composite $$ BG \stackrel{B\pi_G} \longrightarrow B\Z/p\times
B\Z/p \stackrel{yy'} \longrightarrow B^2\Z/p $$ is null-homotopic.
Therefore $B\pi_G$ lifts to $B\phi\colon BG\rightarrow BP(p,3)$,
giving the commutative diagram $$ \xymatrix{
 & & & & BP(p,3) \ar[d]^{B\pi} \\
BP(p,3)\ar[rr]^{\operatorname{res}}& & BG \ar@{.>}[rru]^{B\phi}
\ar[rr]^{B\pi_G} & & B\Z/p\times B\Z/p } $$ which implies that
$B\phi^*(y)=\psi(y)$ and $B\phi^*(y')=\psi(y')$, and
$$(\operatorname{res}\circ\,B\phi)(y)=(\operatorname{res}^*\circ\,\psi)(y)=y\text{
and }
(\operatorname{res}\circ\,B\phi)^*(y')=(\operatorname{res}^*\circ\,\psi)(y)=y'.$$
Now, a standard group theoretical argument (or applying Theorems
\ref{isoP(3,3)} and \ref{isoP(p,3)}) shows that $\phi|_{P(p,3)}$
is an automorphism of $P(p,3)$, that is, $P(p,3)$ is a complement
for $N=\ker\phi\unlhd G$.
\end{proof}

We now proceed with the case $n>3$.

\begin{proposition}\label{prop:P(p,4)}
Let $p$ be an odd prime and $G$ be a finite group such that
$P(p,n)\leq G$.
\begin{enumerate}[\rm (a)]
\item If $n=4$ and there exists $\psi:H^*(BP(p,4))\to H^*(BG)$ as unstable
algebras such that
$(\operatorname{res}\circ\,\psi)|_{H^1_\beta(BP(p,4n))}$ is the
identity, then $P(p,4)$ is a complement for some $N\unlhd G$
\item If $n\geq 5$ and there
exists $\psi:H^*_\beta(BP(p,n))\to H^*_\beta(BG)$ in $\K_\beta$
such that $(\operatorname{res}\circ\,\psi)|_{H^1_\beta(BP(p,n))}$
is the identity, then $P(p,n)$ is a complement for some $N\unlhd
G$
\end{enumerate}
\end{proposition}
\begin{proof}
We now follow along the lines of th proof of Corollary
\ref{cor:P(p,4)}.

If $n=4$ we have $$ 0\rightarrow \Z/p \rightarrow P(p,4)
\stackrel{\pi_1}{\rightarrow} \Z/p\times\Z/p\times\Z/p \rightarrow
0 $$ Let now $B\pi_{1,G}$ be the map $B\pi_{1,G}\colon
BG\rightarrow B\Z/p\times B\Z/p\times B\Z/p$ that classifies the
classes $\psi(u),\psi(y),\psi(y')\in H^1(BG)$.

The composite $$ BG \stackrel{\pi_{1,G}} \longrightarrow B\Z/p
\times B\Z/p\times B\Z/p \stackrel{yy'-\beta(u)} \longrightarrow
B^2\Z/p $$ is null-homotopic because of Remark \ref{remark:tower},
and
$$B\pi_{1,G}^*(yy'-\beta(u))=B\pi_{1,G}^*(y)B\pi_{1,G}^*(y')-B\pi_{1,G}^*(\beta(u))=\psi(y)\psi(y')-\psi(\beta(u))=\psi(yy'-\beta(u))=\psi(0)=0.$$

Therefore $B\pi_{1,G}$ lifts to $B\phi\colon BG\rightarrow
BP(p,4)$, giving the commutative diagram $$\xymatrix{
 & &  & &
BP(p,4) \ar[d]^{B\pi_1} \\BP(p,3)\ar[rr]^{\operatorname{res}}& &
BG \ar@{.>}[rru]^{B\phi} \ar[rr]^-{B\pi_{1,G}} & & B\Z/p\times
B\Z/p\times B\Z/p.} $$ This implies that $B\phi^*(y)=\psi(y)$,
$B\phi^*(y')=\psi(y')$, and $B\phi^*(u)=\psi(u)$, and
$$(\operatorname{res}\circ\,B\phi)(y)=(\operatorname{res}^*\circ\,\psi)(y)=y,$$
$$(\operatorname{res}\circ\,B\phi)^*(y')=(\operatorname{res}^*\circ\,\psi)(y)=y',\text{
and }$$
$$(\operatorname{res}\circ\,B\phi)(u)=(\operatorname{res}^*\circ\,\psi)(u)=u.$$
Now, a standard group theoretical argument (or applying Theorem
\ref{isoP(p,n)}) shows that $\phi|_{P(p,4)}$ is an automorphism of
$P(p,4)$, that is, $P(p,4)$ is a complement for $N=\ker\phi\unlhd
G$.

If $n>4$, we must consider before the extensions: $$ 0\rightarrow
\Z/p \longrightarrow \Z/p^{i}\times\Z/p\times\Z/p
\stackrel{\pi_{i-1}}{\longrightarrow}
\Z/p^{i-1}\times\Z/p\times\Z/p \rightarrow 0 , $$ classified by
$\beta_{i-1}(u) \in H^*_\beta(\Z/p^{i-1}\times\Z/p\times\Z/p)$.

Let $B\pi_{1,G}\colon BG\rightarrow B\Z/p\times B\Z/p\times\Z/p$
be the map that classifies the classes
$\psi(u),\psi(y),\psi(y')\in H^1(BG)$. As $\beta_1(\psi(u))=0 \in
H^*_\beta(BG)$ the map $B\pi_{1,G}$ extends to a map $B\pi_{2,G}$
having the following commutative diagram: $$ \xymatrix{
  & & B\Z/p^2\times B\Z/p \times B\Z/p \ar[d]^{B\pi_1} \\
BG \ar@{.>}[rru]^{B\pi_{2,G}} \ar[rr]^-{B\pi_{1,G}} & &
B\Z/p\times B\Z/p\times B\Z/p
} $$ Using the same argument we can proceed extending the map till
$B\pi_{n-3} \colon B\Z/p^{n-3} \times B\Z/p \times B\Z/p$. To do
the last step we use again Remark \ref{remark:tower} and that
$\beta_{n-3}(\psi(u))=\psi(y)\psi(y') \in H^*_\beta(BG)$, and
obtain a map $B\phi$ which gives the commutative diagram: $$
\xymatrix{
  & & BP(p,n) \ar[d]^{B\pi_{n-2} \circ \cdots \circ B\pi_1} \\
BG \ar@{.>}[rru]^{B\phi} \ar[rr]^-{B\pi_{1,G}} & & B\Z/p\times
B\Z/p\times B\Z/p} $$ which implies again that
$B\phi^*(y)=\psi(y)$, $B\phi^*(y')=\psi(y')$, and
$B\phi^*(u)=\psi(u)$, and
$$(\operatorname{res}\circ\,B\phi)(y)=(\operatorname{res}^*\circ\,\psi)(y)=y,$$
$$(\operatorname{res}\circ\,B\phi)^*(y')=(\operatorname{res}^*\circ\,\psi)(y)=y',\text{
and }$$
$$(\operatorname{res}\circ\,B\phi)(u)=(\operatorname{res}^*\circ\,\psi)(u)=u.$$
Again, a standard group theoretical argument, or applying Theorem
\ref{isoP(p,n)}, we obtain that $\phi|_{P(p,n)}$ is an
automorphism of $P(p,n)$, that is, $P(p,n)$ is a complement for
$N=\ker\phi\unlhd G$.
\end{proof}


\end{document}